\input amstex
\input amsppt.sty
\magnification=\magstep1
\hsize=33truecc
\vsize=22.2truecm
\baselineskip=16truept
\NoBlackBoxes
\nologo
\pageno=1
\topmatter
\TagsOnRight

\def\N{\Bbb N}
\def\Z{\Bbb Z}

\def\l{\left}
\def\r{\right}
\def\b{\bigg}

\def\({\b(}
\def\[{\b[}
\def\){\b)}
\def\]{\b]}

\def\t{\text}
\def\f{\frac}
\def\ord{\roman{ord}}
\def\mo{\roman{mod}}

\def\se {\subseteq}

\def\sm{\setminus}

\def\bi{\binom}
\def\eq{\equiv}

\def\ls{\leqslant}
\def\gs{\geqslant}

\def\ve{\varepsilon}

\def\Proof{\noindent{\it Proof}}
\def\Remark{\noindent{\it Remark}}

\def\Ack{\noindent {\bf Acknowledgment}}
\hbox{J. Number Theory 162(2016), 190--211.}
\medskip
\title A result similar to Lagrange's theorem\endtitle
\author Zhi-Wei Sun \endauthor
\affil Department of Mathematics, Nanjing University
     \\Nanjing 210093, People's Republic of China
    \\  zwsun\@nju.edu.cn
    \\ {\tt http://math.nju.edu.cn/$\sim$zwsun}
 \endaffil
\abstract Generalized octagonal numbers are those $p_8(x)=x(3x-2)$ with $x\in\Z$.
In this paper we show that every positive integer can be written as the sum of four generalized octagonal numbers one of which is odd.
This result is similar to Lagrange's theorem on sums of four squares. Moreover, for $35$ triples $(b,c,d)$ with $1\ls b\ls c\ls d$ (including $(2,3,4)$ and $(2,4,8)$),
we prove that any nonnegative integer can be exprssed as $p_8(w)+bp_8(x)+cp_8(y)+dp_8(z)$ with $w,x,y,z\in\Z$. We also pose several conjectures for further research.
\endabstract
\thanks 2010 {\it Mathematics Subject Classification}.
Primary 11E25; Secondary 11B75, 11D85, 11E20.
\newline \indent Supported by the National Natural Science
Foundation of China (grant 11571162).
\endthanks
\endtopmatter
\document

\heading{1. Introduction}\endheading

Polygonal numbers are nonnegative integers constructed geometrically from the regular polygons.
For $m=3,4,\ldots$, those {\it $m$-gonal numbers} (or {\it polygonal numbers of order $m$}) are given
by
$$p_m(n):=(m-2)\bi n2+n=\f{(m-2)n^2-(m-4)n}2\ (n=0,1,2,\ldots),$$
and those $p_m(x)$ with $x\in\Z$ are called the {\it generalized $m$-gonal numbers}. Note that
$$\gather p_3(x)=\f{x(x+1)}2,\ p_4(x)=x^2,\ p_5(x)=\f{x(3x-1)}2,
\\ p_6(x)=x(2x-1),\ p_7(x)=\f{x(5x-3)}2,\ p_8(x)=x(3x-2).
\endgather$$
Fermat's claim that each $n\in\N=\{0,1,2,\ldots\}$  can be written as the sum of $m$ polygonal numbers of order $m$
was proved by Lagrange in the case $m=4$, Gauss in the case $m=3$, and Cauchy in the case $m\gs5$. (See, e.g., [N, pp.\, 3-35] and [MW, pp.\,54-57].)

It is easy to see that generalized hexagonal numbers coincide with triangular numbers. In 1994, R. K. Guy [Gu] observed that
each nonnegative integer can be written as the sum of three generalized pentagonal numbers.
By a theorem of Legendre (cf. [N, p.\,33]), for each positive odd number $m\gs5$, any integer $n\gs28(m-2)^3$ can be expressed as the sum of four $m$-gonal numbers;
in particular, any integer $n\gs3500$ is the sum of four heptagonal numbers. Via a computer we find that every $n=0,1,\ldots,3500$ can be written as the sum of four generalized heptagonal numbers. So we have
$$\{p_7(w)+p_7(x)+p_7(y)+p_7(z):\ w,x,y,z\in\Z\}=\N.\tag1.1$$
For any integer $m>8$, clearly $5$ cannot be written as the sum of four generalized $m$-gonal numbers.

Octagonal numbers are those $p_8(n)=n(3n-2)$ with $n\in\N$, and generalized octagonal numbers are
$p_8(x)=x(3x-2)$ with $x\in\Z$. The sequence of generalized octagonal numbers appears as [SD, A001082]. Here is the list of generalized octagonal numbers up to $120$:
$$0,\ 1,\ 5,\, 8,\, 16,\ 21,\ 33,\ 40,\ 56,\ 65,\ 85,\ 96,\ 120.$$

In this paper we establish the following new theorem which is quite similar to Lagrange's theorem on sums of four squares.

\proclaim{Theorem 1.1} Let $n$ be any positive integer. Then $n$ can be written as the sum of four generalized octagonal numbers one of which is odd, i.e.,
there are $w,x,y,z\in\Z$ not all even such that
$$n=w(3w-2)+x(3x-2)+y(3y-2)+z(3z-2).\tag1.2$$
\endproclaim

For $n\in\Z^+=\{1,2,3,\ldots\}$, we let $r(n)$ denote the number of ways to write $n$ as the sum of four unordered generalized octagonal numbers, and define
$s(n)$ to be the number of ways to write $n$ as the sum of four unordered generalized octagonal numbers not all even. Clearly, $r(n)\gs s(n)$ for all $n\in\Z^+$.
\medskip

{\it Example}\ 1.1. We have $r(n)=s(n)=1$ for $n=1,3,5,9,13$; in fact,
$$1=0+0+0+1,\ 3=0+1+1+1,\ 5=0+0+0+5,\ 9=0+0+1+8,\ 13=0+0+5+8.$$
Also, $r(n)>s(n)=1$ for $n=8,16,24,40,56$; in fact,
$$\gather 8=1+1+1+5=0+0+0+8,
\\ 16=1+5+5+5=0+0+0+16=0+0+8+8,
\\\ 24=1+1+1+21=0+0+8+16=0+8+8+8,
\\40=1+1+5+33=0+0+0+40=0+8+16+16=8+8+8+16,
\endgather$$
and
$$56=1+1+21+33=0+0+16+40=0+8+8+40=8+8+16+24=8+16+16+16.$$

We will prove Theorem 1.1 in the next section and study the sets $\{n\in\Z^+:\ r(n)=1\}$ and $\{n\in\Z^+:\ s(n)=1\}$ in Section 3.

Let $f_1(x),\ldots,f_k(x)$ be integer-valued polynomials, and let $a_1,\ldots,a_k$ be positive integers. If any $n\in\N$ can be written as
$$a_1f_1(x_1)+\cdots+a_kf_k(x_k)\ \ \t{with}\ x_1,\ldots,x_k\in\N,$$
then we call $a_1f_1+\cdots+a_kf_k$ a {\it universal sum over $\N$}.
Similarly, if any $n\in\N$ can be written as
$$a_1f_1(x_1)+\cdots+a_kf_k(x_k)\ \ \t{with}\ x_1,\ldots,x_k\in\Z,$$
then we call $a_1f_1+\cdots+a_kf_k$ a {\it universal sum over $\Z$}.
In 1862 Liouville (cf. [D99, p.\,23]) determined all those universal sums $ap_3+bp_3+cp_3$ with $a,b,c\in\Z^+$.
In 1917 Ramanujan [R] listed all the 54 universal sums $ap_4+bp_4+cp_4+dp_4$ with $1\ls a\ls b\ls c\ls d$, and
the list was later confirmed by Dickson [D27].
The author [S15] systematically investigated universal sums $ap_i+bp_j+cp_k$ (over $\N$ or $\Z$)
with $a,b,c\in\Z^+$ and $i,j,k\in\{3,4,5,\ldots\}$; for example, he proved that $p_3+4p_4+p_5$ is universal over $\Z$ and conjectured that it is even universal over $\N$.

Since
 $$p_8(2x+1)=(2x+1)(6x+1)=4p_8(-x)+1,\tag1.3$$
Theorem 1.1 indicates that both $p_8+p_8+p_8+p_8$ and $p_8+p_8+p_8+4p_8$ are universal over $\Z$.
Motivated by this, we aim to find all universal sums $ap_8+bp_8+cp_8+dp_8$ over $\Z$, where $a,b,c,d$ are positive integers.

\proclaim{Theorem 1.2} Let $a,b,c,d\in\Z^+$ with $a\ls b\ls c\ls d$.
Suppose that $ap_8+bp_8+cp_8+dp_8$ is universal over $\Z$. Then we must have $a=1$, and $(b,c,d)$ is
among the following $40$ triples:
$$\aligned&(1,1,1),\,(1,1,2),\,(1,1,3),\,(1,1,4),\,(1,2,2),\,(1,2,3),\,(1,2,4),\,(1,2,5),
\\&(1,2,6),\,(1,2,7),\,(1,2,8),\,(1,2,9),\,(1,2,10),\,(1,2,11),\,(1,2,12),\,(1,2,13),
\\&(1,3,3),\,(1,3,5),\,(1,3,6),\,(2,2,2),\,(2,2,3),\,(2,2,4),\,(2,2,5),\,(2,2,6),
\\&(2,3,4),\,(2,3,5),\,(2,3,6),\,(2,3,7),\,(2,3,8),\,(2,3,9),\,(2,4,4),\,(2,4,5),
\\&(2,4,6),\,(2,4,7),\,(2,4,8),\,(2,4,9),\,(2,4,10),\,(2,4,11),\,(2,4,12),\,(2,4,13).
\endaligned$$
\endproclaim

\proclaim{Theorem 1.3}  $p_8+bp_8+cp_8+dp_8$ is universal over $\Z$ for any $(b,c,d)$ among the $33$ triples
$$\aligned&(1,2,2),\,(1,2,8),\,(2,2,4),\,(2,4,8),\,(2,2,2),\,(2,4,4),
\\&(1,1,2),\,(1,2,3),\,(1,2,5),\,(1,2,7),\,(1,2,9),\,(1,2,11),\,(1,2,13),
\\&(1,2,4),\,(2,3,4),\,(2,4,5),\,(2,4,7),\,(2,4,9),\,(2,4,11),\,(2,4,13),
\\&(1,1,3),\,(2,2,3),\,(2,2,6),\,(2,3,8),\,(1,2,6),\,(1,2,10),\,(1,2,12),
\\&(2,4,6),(2,4,10),\,(2,4,12),\,(2,2,5),\,\,(2,3,5),\,(1,3,5).
\endaligned$$
\endproclaim

Theorems 1.2 and 1.3 will be proved in Section 4. Below is our related conjecture.

\proclaim{Conjecture 1.1} {\rm (i)} $p_8+bp_8+cp_8+dp_8$ is universal over $\Z$ if $(b,c,d)$ is among the five triples
$$\aligned &(1,3,3),\ (1,3,6),\ (2,3,6),\ (2,3,7),\ (2,3,9).
\endaligned\tag1.4$$

{\rm (ii)} If $(b,c,d)$ is among the five triples
$$(1,1,2),\ (1,2,3),\ (1,2,5),\ (1,2,11),\ (2,3,4),$$
then any $n\in\N$ can be written as $p_8(w)+bp_8(x)+cp_8(y)+dp_8(z)$ with $w\in\Z$ and $x,y,z\in\N$.
\endproclaim
\Remark\ 1.1. Part (i) of Conjecture 1.1 is a supplement to Theorems 1.1-1.3. In fact, if we remove from
the 40 triples in Theorem 1.2, the two triples $(1,1,1)$, $(1,1,4)$ treated in Theorem 1.1 and the 33 triples in Theorem 1.3, then we get the remaining 5 triples listed in (1.4).
\medskip

Here is another conjecture.

\proclaim{Conjecture 1.2} We have
$$\align\{p_8(x)+p_8(y)+3p_8(z):\ x,y,z\in\Z\}=&\N\sm\{7,14,18,91\},\tag1.5
\\\{p_8(x)+p_8(y)+6p_8(z):\ x,y,z\in\Z\}=&\N\sm\{3,4,18,20,25,108,298\},
\\\{p_8(x)+3p_8(y)+3p_8(z):\ x,y,z\in\Z\}\supseteq& \{n\in\N:\ n>3265\},
\\\{p_8(x)+3p_8(y)+7p_8(z):\ x,y,z\in\Z\}\supseteq& \{n\in\N:\ n>1774\},
\\\{p_8(x)+3p_8(y)+9p_8(z):\ x,y,z\in\Z\}\supseteq& \{n\in\N:\ n>446\}.
\endalign$$
\endproclaim
\Remark\ 1.2. For any $n\in\N$, it is easy to see that $n=p_8(u)+p_8(v)+3p_8(w)$ for some $u,v,w\in\Z$ if and only if
$3n+5=x^2+y^2+3z^2$ for some $x,y,z\in\Z$ with $3\nmid z$. Thus, (1.5) holds if and only if for any $n\in\Z^+$ with $n\not=8,15,19,92$
we can write $3n+2=x^2+y^2+3z^2$ with $x,y,z\in\Z$ and $3\nmid z$. We also conjecture that for any $n\in\Z^+$ with $n\not=3,10,11,55,150$
we can write $3n+1=x^2+y^2+3z^2$ with $x,y,z\in\Z$ and $3\nmid z$.
\medskip
Inspired by Theorems 1.1-1.3 and Conjectures 1.1-1.2, we are going to pose in Section 5 some similar conjectures involving $p_m(x)$ with $m\in\{5,6,7\}$.

\heading{2. Proof of Theorem 1.1}\endheading

\proclaim{Lemma 2.1} Any integer $n>4$ can be written as the sum of four squares one of which is even and two of which are nonzero.
\endproclaim
\Proof. It is well-known that
$$r_4(m)=8\sum\Sb d\mid m\\4\nmid d\endSb d\quad\t{for all}\ m=1,2,3,\ldots,\tag2.1$$
where
$$r_4(m):=|\{(w,x,y,z)\in\Z^4:\ w^2+x^2+y^2+z^2=n\}|.$$
(See, e.g., [B, p.\,59].) If $m>1$ is an integer whose smallest prime divisor is $p$, then
$$r_4(m)\gs 8(1+p)>2^4$$
and hence $m$ can be written as the sum of four squares (at least) two of which are nonzero.

By the above, we can write any integer $n>4$ as the sum of four squares two of which are nonzero. If all the four squares are odd, then
$n\eq4\pmod 8$ and we can write $n/4>1$ in the form $w^2+x^2+y^2+z^2$ with $w,x,y,z\in\Z$ and $wx\not=0$, hence
$n=(2w)^2+(2x)^2+(2y)^2+(2z)^2$ with $2w\not=0$ and $2x\not=0$. This completes the proof. \qed

\medskip
\Remark\ 2.1. As first stated by J. Liouville [L] in 1861 and proved by T. Pepin [P] in 1890,
for $n=2^an_0$ with $a\in\N$ and $n_0\in\{1,3,5,\ldots\}$,
we have
$$\align&|\{(w,x,y,z)\in\Z^4:\ w^2+x^2+y^2+4z^2=n\}|
\\=&\cases 2(2+(-1)^{(n_0-1)/2})\sigma(n_0)&\t{if}\ a=0,
\\12\sigma(n_0)&\t{if}\ a=1,
\\8\sigma(n_0)&\t{if}\ a=2,
\\24\sigma(n_0)&\t{if}\ a>2,
\endcases
\endalign$$
where $\sigma(n_0)$ is the sum of all positive divisors of $n_0$.  See [AALW, Theorem 1.7] for this known result and some other similar ones.

\proclaim{Lemma 2.2} {\rm (i)} Suppose that $x,y,z\in\Z$ are not all divisible by $3$. Then there are $\bar x,\bar y,\bar z\in\Z$ not all divisible by $3$ such that
$$\bar x\eq x\ (\mo\ 2),\ \bar y\eq y\ (\mo\ 2),\ \bar z\eq z\ (\mo\ 2),\ \t{and}\ 9(x^2+y^2+z^2)=\bar x^2+\bar y^2+\bar z^2.$$

{\rm (ii)} Suppose that $x,y,z$ are integers with $x^2+y^2+z^2$ a positive multiple of $3$. Then $x^2+y^2+z^2=\bar x^2+\bar y^2+\bar z^2$ for some $\bar x,\bar y,\bar z\in\Z$ with
$$\bar x\eq x\ (\mo\ 2),\ \bar y\eq y\ (\mo\ 2),\ \bar z\eq z\ (\mo\ 2),\ \t{and}\ 3\nmid \bar x\bar y\bar z.\tag2.2$$
\endproclaim
\Proof. (i) As $x,y,z$ are not all divisible by $3$, there are $x'\in\{\pm x\}$, $y'\in\{\pm y\}$ and $z'\in\{\pm z\}$ such that
$x'+y'+z'\not\eq0\pmod 3$.
Let $\bar x=x'-2y'-2z'$, $\bar y=y'-2x'-2z'$ and $\bar z=z'-2x'-2y'$. It is easy to verify the identity
$$9((x')^2+(y')^2+(z')^2)=\bar x^2+\bar y^2+\bar z^2$$
which is a special case of R\'ealis' identity (cf. [D99, p.\,266]). Clearly,
$$\bar x\eq x\ (\mo\ 2),\ \bar y\eq y\ (\mo\ 2),\ \bar z\eq z\ (\mo\ 2),$$
and
$$\bar x\eq \bar y\eq \bar z\eq x'+y'+z'\not\eq0\pmod 3.$$
This proves part (i).

(ii) Let $a\in\N$ be the $3$-adic order of $\gcd(x,y,z)$, and write $x=3^ax_0$, $y=3^ay_0$ and $z=3^az_0$, where $x_0,y_0,z_0$ are integers not all divisible by $3$.
Note that $x^2+y^2+z^2=9^a(x_0^2+y_0^2+z_0^2)$. Applying part (i) again and again, we finally get that
$$9^a(x_0^2+y_0^2+z_0^2)=\bar x^2+\bar y^2+\bar z^2$$
for some $\bar x,\bar y,\bar z\in\Z$ not all divisible by $3$ with
$$\bar x\eq x_0\eq x\pmod2,\ \bar y\eq y_0\eq y\pmod2,\ \t{and}\ \bar z\eq z_0\eq z\pmod3.$$
As $\bar x^2+\bar y^2+\bar z^2=x^2+y^2+z^2\eq0\pmod 3$, we must have $3\nmid \bar x\bar y\bar z$. This concludes the proof of part (ii). \qed

\medskip
\Remark\ 2.2. R\'ealis' identity discovered in 1878 is as follows (cf. [D99, p.\,266]):
$$\align &(a^2+b^2+c^2)^2(x^2+y^2+z^2)=((b^2+c^2-a^2)x-2a(by+cz))^2
\\&+((a^2-b^2+c^2)y-2b(ax+cz))^2+((a^2+b^2-c^2)z-2c(ax+by))^2.
\endalign$$

\medskip
\noindent{\it Proof of Theorem} 1.1. For $w,x,y,z\in\Z$, we clearly have
$$\align &n=w(3w-2)+x(3x-2)+y(3y-2)+z(3z-2)
\\\iff &3n+4=(3w-1)^2+(3x-1)^2+(3y-1)^2+(3z-1)^2.
\endalign$$
If an integer $m$ is not divisible by $3$, then $m$ or $-m$ can be written as $3x-1$ with $x\in\Z$.
Also, $(3(1-2x)-1)^2=4(3x-1)^2$ for any $x\in\Z$.
Thus, it suffices to show that $3n+4$ can be written as the sum of four squares none of which is divisible by $3$ and one of which is even.

By Lemma 2.1, we may write $3n+4$ as $w^2+x^2+y^2+z^2$, where $w,x,y,z$ are integers one of which is even and two of which are nonzero.
Clearly, $w,x,y,z$ cannot be all divisible by $3$. Without loss of generality, we suppose that $3\nmid w$. Note that $x,y,z$ are not all zero and $x^2+y^2+z^2\eq 4-w^2\eq0\pmod3$.
By Lemma 2.2(ii), $x^2+y^2+z^2=\bar x^2+\bar y^2+\bar z^2$ for some $\bar x,\bar y,\bar z\in\Z$ satisfying (2.2).
Clearly $2\mid w\bar x\bar y\bar z$ since $2\mid wxyz$. Note that $3n+4=w^2+\bar x^2+\bar y^2+\bar z^2$ and $3\nmid w\bar x\bar y\bar z$.
This concludes our proof. \qed

\heading{3. On the sets $\{n\in\Z^+:\ r(n)=1\}$ and $\{n\in\Z^+:\ s(n)=1\}$}\endheading

In view of the first paragraph in the proof of Theorem 1.1 given in the last section, for any positive integer $n$, $r(n)=1$ if and only if
$3n+4$ can be written uniquely as the sum of four squares not divisible by $3$. Also, $s(n)=1$ if and only if
$3n+4$ can be written uniquely as the sum of four squares all coprime to $3$ but not all odd.

{\it Example}\ 3.1. $r((2^{2k}-4)/3)=s((2^{2k}-4)/3=1$ for any integer $k>1$. This is because
$$\align 2^{2k}=&(\pm2^{k-1})^2+(\pm2^{k-1})^2+(\pm2^{k-1})^2+(\pm2^{k-1})^2
=(\pm2^k)^2+0+0+0
\\=&0+(\pm2^k)^2+0+0=0+0+(\pm2^k)^2+0=0+0+0+(\pm2^k)^2
\endalign$$
and $r_4(2^{2k})=8(1+2)=16+8$ by (2.1).
\medskip

{\it Example}\ 3.2. $r((2^{2n+1}5-4)/3)=s((2^{2n+1}5-4)/3)=1$ for any $n\in\N$. In fact, as
$$\align 2^{2n+1}5=&(\pm2^{n})^2(1^2+1^2+2^2+2^2)=(\pm2^n)^2(0^2+0^2+1^2+3^2)
\endalign$$
and $r_4(2^{2n+1}5)=8(1+2+5+10)=\bi{4}22^4+2\bi 422^2$ by (2.1), the only ways to write
$2^{2n+1}5$ in the form $w^2+x^2+y^2+z^2$ with $w,x,y,z\in\N$ and $w\ls x\ls y\ls z$ are
$$2^{2n+1}5=(2^n)^2+(2^n)^2+(2^{n+1})^2+(2^{n+1})^2=0^2+0^2+(2^n)^2+(2^n3)^2.$$

{\it Example}\ 3.3. $r((2^{2n+1}11-4)/3)=s((2^{2n+1}11-4)/3)=1$ for any $n\in\N$. In fact, as
$$\align 2^{2n+1}11=&(\pm2^n)^2(1^2+1^2+2^2+4^2)=(\pm2^{n})^2(0^2+2^2+3^2+3^2)
\endalign$$
and $r_4(2^{2n+1}11)=8(1+2+11+22)=2\bi 422^4+2\bi 422^3$ by (2.1), the only ways to write
$2^{2n+1}11$ in the form $w^2+x^2+y^2+z^2$ with $w,x,y,z\in\N$ and $w\ls x\ls y\ls z$ are
$$2^{2n+1}11=(2^n)^2+(2^n)^2+(2^{n+1})^2+(2^{n+2})^2=0^2+(2^{n+1})^2+(2^{n}3)^2+(2^{n}3)^2.$$

{\it Example}\ 3.4. $r((2^{2n+1}23-4)/3)=s((2^{2n+1}23-4)/3)=1$ for any $n\in\N$. In fact, as
$$\align 2^{2n+1}23=&(\pm2^n)^2(1^2+2^2+4^2+5^2)=(\pm2^{n})^2(0^2+1^2+3^2+6^2)
\endalign$$
and $r_4(2^{2n+1}23)=8(1+2+23+46)=2^4\times 4!+2^3\times 4!$ by (2.1), the only ways to write
$2^{2n+1}23$ in the form $w^2+x^2+y^2+z^2$ with $w,x,y,z\in\N$ and $w\ls x\ls y\ls z$ are
$$2^{2n+1}23=(2^n)^2+(2^{n+1})^2+(2^{n+2})^2+(2^{n}5)^2=0^2+(2^{n})^2+(2^{n}3)^2+(2^{n+1}3)^2.$$

In view of Example 1.1 and Examples 3.1-3.4, we propose the following conjecture based on our computation.

\proclaim{Conjecture 3.1} Let $n$ be a positive integer.
If $r(n)=1$, then we must have
$$3n+4\in\{7,\,13,\,19,\,31,\,43\}\cup E$$
where
$$E:=\{2^{2k}:\ k=2,3,\ldots\}\cup\bigcup_{n\in\N}\l\{2^{2n+1}5,\ 2^{2n+1}11,\ 2^{2n+1}23\r\}.$$
If $s(n)=1$, then we must have
$$3n+4\in\{7,\,13,\,19,\,31,\,43,\,4\times7,\,4\times13,\,4\times 19,\,4\times 31,\,4\times43\}\cup E.$$
\endproclaim

\heading{4. Proofs of Theorems 1.2 and 1.3}\endheading

Recall that
$$\{p_8(x):\ x\in\Z\}=\{n(3n\pm2):\ n=0,1,2,\ldots\}=\{0,1,5,8,16,21,33,\ldots\}.$$

\medskip
\noindent{\it Proof of Theorem} 1.2. For convenience, we define
$$Z(a,b,c,d):=\{ap_8(w)+bp_8(x)+cp_8(y)+dp_8(z):\ w,x,y,z\in\Z\}.$$

As $1\in Z(a,b,c,d)$ and $a\ls b\ls c\ls d$, we must have $a=1$. Note that $b\ls 2$ since $2\in Z(1,b,c,d)$.

{\it Case} 1. $b=1$.

By $3\in Z(1,1,c,d)$, we must have $c\ls 3$.

If $c=1$, then $d\ls 4$ by $4\in Z(1,1,1,d)$.

When $c=2$, we must have $d\ls 14$ by $14\in Z(1,1,2,d)$. Note that $Z(1,1,2,14)$ does not contain $60$.

If $c=3$, then $d\ls 7$ by $7\in Z(1,1,3,d)$. Note that
$$18\not\in Z(1,1,3,4)\ \t{and}\ 14\not\in Z(1,1,3,7).$$

{\it Case} 2. $b=2$.

As $p_8(x)+2p_8(y)\not=4$ for any $x,y\in\Z$, we have $c\ls 4$ by $4\in Z(1,2,c,d)$.

If $c=2$, then $d\ls 6$ by $6\in Z(1,2,2,d)$.

When $c=3$, we have $d\ls 9$ by $9\in Z(1,2,3,d)$.
Note that $12\not\in Z(1,2,3,3)$.

If $c=4$, then $d\ls 13$ by $13\in Z(1,2,4,d)$.

Combining the above, we obtain the desired result. \qed

\proclaim{Lemma 4.1} A positive integer $n$ can be written as the sum of four nonzero squares, if and only if it does not belong to the set
$$\{1,\,3,\,5,\,9,\,11,\,17,\,29,\,41\}\cup\bigcup_{k\in\N}\{2\times 4^k,\,6\times 4^k,\,14\times 4^k\}.$$
\endproclaim
\Remark\ 4.1. This is a known result, see, e.g., [G, pp.\,74--75].

\proclaim{Lemma 4.2} Let $w=x^2+my^2$ be a positive integer with $m\in\{2,5,8\}$ and $x,y\in\Z$. Then $w=u^2+mv^2$ for some integers $u$ and $v$ not all divisible by $3$.
\endproclaim
\Remark\ 4.2. This is [S15, Lemma 2.1].
\medskip

The famous Gauss-Legendre theorem on sums of three squares (cf. [N, p.\,23]) asserts that
$$\{x^2+y^2+z^2:\ x,y,z\in\Z\}=\N\sm\{4^k(8l+7):\ k,l\in\N\}.\tag4.1$$

\proclaim{Lemma 4.3} {\rm (i)} Any positive odd integer can be expressed as $x^2+y^2+2z^2$ with $x,y,z\in\Z$.

{\rm (ii)} For any positive integer $n$, we can write $6n+1$ as $x^2+y^2+2z^2$, where $x,y,z$ are integers with $2\mid xy$ and $3\nmid xyz$.
\endproclaim
\Proof. (i) Part (i) was first observed by Euler (cf. [D99, p.\,260]). In fact, by (4.1),
for any $n\in\N$ we can write $4n+2$ as the sum of three squares.
Thus, there are $x,y,z\in\Z$ such that
$$4n+2=(2x+1)^2+(2y+1)^2+(2z)^2=2(x+y+1)^2+2(x-y)^2+4z^2$$
and hence $2n+1=(x+y+1)^2+(x-y)^2+2z^2$.

(ii) Let $n\in\Z^+$. If $6n+1=m^2$ for some (odd) integer $m>1$, then by S. Cooper and H. Y. Lam [CL] we have
$$\align&|\{(x,y,z)\in\Z^3:\ x^2+y^2+2z^2=m^2\}|
\\=&4\prod_{p>2}\f{p^{\ord_p(m)+1}-1-(\f{-2}p)(p^{\ord_p(m)}-1)}{p-1}\gs 4\prod_{p>2}p^{\ord_p(m)}>4,
\endalign$$
where $\ord_p(m)$ stands for the order of $m$ at the prime $p$, and $(\f{\cdot}p)$ denotes the Legendre symbol.
In view of this and part (i), we can always write $6n+1=x^2+y^2+2z^2$ with $x,y,z\in\Z$ and $x^2,y^2\not=6n+1$. As $2z^2\not\eq1\pmod3$,
one of $x$ and $y$, say $x$, is not divisible by $3$. Since $y^2+2z^2$ is a positive multiple of $3$, by Lemma 4.2 we can
write $y^2+2z^2=\bar y^2+2\bar z^2$ with $\bar y,\bar z\in\Z$ and $3\nmid \bar y\bar z$. Thus
$6n+1=x^2+\bar y^2+2\bar z^2$ with $3\nmid x\bar y\bar z$. Clearly $x\not\eq\bar y\pmod2$ and hence $2\mid x\bar y$. This proves part (ii). \qed

\proclaim{Lemma 4.4} Let $n\in\N$ and $r\in\{1,3,5,7\}$.
Let $a,b,c,d$ be integers with
$$a\eq1\ (\mo\ 2),\ b\eq2\ (\mo\ 4),\ c\eq0\ (\mo\ 4)\ \t{and}\ d\eq r\ (\mo\ 4).$$

{\rm (i)} If $d\not\eq r\pmod 8$, then for some $w\in\{a,b,c\}$ we have $n+dw^2\not=4^k(8m+r)$ for all $k\in\N$ and $m\in\Z$.

{\rm (ii)} We have $n-dw^2\not\in S$ for some $w\in\{a,b,c\}$, where
$$S:=\{8q-d:\ q\in\Z\}\cup\{4^k(8l+r):\ k,l\in\N\}.\tag4.2$$
\endproclaim
\Proof. Clearly, $a^2,b^2,c^2$ are congruent to $1,4,0$ modulo $8$ respectively.
Thus $da^2,db^2,dc^2$ are pairwise incongruent modulo $8$. Note that
$$db^2\eq dc^2\not\eq da^2\pmod 4.$$
For any $m\in\Z$, obviously $4(8m+r)\eq 4r\eq4\pmod 8$
and $4^k(8m+r)\eq0\pmod 8$ for $k=2,3,\ldots$.

(i) Now assume that $d\not\eq r\pmod 8$ and
$$\{n+da^2,n+db^2,n+dc^2\}\se \{4^k(8m+r):\,k\in\N,\ m\in\Z\}.$$
We want to deduce a contraction. By the above analysis, we must have $n+da^2\eq r\pmod 8$. Hence $n+dc^2\eq r+d(c^2-a^2)\eq r+d(0-1)\eq 4\pmod 8$ and thus
$n+dc^2=4(8q+r)$ for some $q\in\Z$. Note that $d+r\eq2r\eq2\pmod 4$ and hence
$$n+db^2=4(8q+r)+d(b^2-c^2)\eq 4r+d(4-0)\eq8\pmod {16},$$
which contradicts that $n+db^2\in\{4^k(8m+r):\ k\in\N,\ m\in\Z\}$. This proves part (i).

(ii) Now we come to show part (ii). Suppose that $\{n-da^2,n-db^2,n-dc^2\}\se S.$
As $n-db^2-(n-dc^2)\eq 4\pmod 8$ and $-d\not\eq r\pmod 4$, we must have $n-db^2,n-dc^2\in\{4^k(8l+r): k\in\Z^+,\ l\in\N\}$.
Hence $n-da^2$ is congruent to $r$ or $-d$ modulo $8$. If $n-da^2\eq r\pmod 8$, then
$$n-db^2\eq r+d(a^2-b^2)\eq r+d\eq 2r\eq 2\pmod 4$$
which contradicts that $n-db^2\in S$. If $n-da^2\eq -d\pmod 8$, then
$$n-db^2\eq -d+d(a^2-b^2)\eq -d+d(1-4)=-4d\eq4\pmod 8$$
and hence $n-db^2=4(8m+r)$ for some $m\in\N$, therefore
$$\align n-dc^2=&4(8m+r)+d(b^2-c^2)
\\\eq& 4(8m+r)+d(4-0)=4(8m+d+r)\pmod {16},
\endalign$$
which contradicts that $n-dc^2\in S$ since $d+r\eq 2r\eq2\pmod 4$.
This proves part (ii). \qed

\medskip
\noindent{\it Proof of Theorem 1.3}. For $b,c,d\in\Z^+$ and $n\in\N$, clearly
$$\align &n=p_8(w)+bp_8(x)+cp_8(y)+dp_8(z)
\\\iff& 3n+b+c+d+1=(3w-1)^2+b(3x-1)^2+c(3y-1)^2+d(3z-1)^2.
\endalign$$
For any integer $m\not\eq0\pmod3$, either $m$ or $-m$ can be written as $3x-1$ with $x\in\Z$.
So, $p_8+bp_8+cp_8+dp_8$ is universal over $\Z$ if and only if for any $n\in\N$ we have
$$3n+b+c+d+1=w^2+bx^2+cy^2+dz^2\ \ \t{for some}\ w,x,y,z\in\Z\ \t{with}\ 3\nmid wxyz.\tag4.3$$
Below we will use this simple fact quite often.

(a) We first prove the universality of $p_8+2p_8+4p_8+8p_8$ over $\Z$.
It suffices to show that for any given $n\in\N$ we have $3n+15=w^2+x^2+2y^2+2z^2$ for some integers $w,x,y,z$ with $2\mid wx$, $2\mid yz$ and $\gcd(wxyz,3)=1$.

 If $3n+15=36$, then the representation $36=2^2+4^2+2\times2^2+2\times 2^2$ suffices. When $3n+15=6\times 4^k$ for some $k\in\Z^+$, the representation
$$3n+15=(2^k)^2+(2^k)^2+2(2^k)^2+2(2^k)^2$$
meets our purpose.

Now we suppose that $3n+15\not=36$ and $3n+15\not=6\times 4^k$ for any $k\in\N$. As $3n+15\eq0\pmod 3$, by Lemma 4.1 we can write $3n+15=w^2+x^2+y^2+z^2$ with $w,x,y,z$ nonzero integers.
When $2\nmid wxyz$, we have $3n+15\eq4\pmod 8$, hence by Lemma 4.1 there are nonzero integers $w_0,x_0,y_0,z_0$ such that
$(3n+15)/4=w_0^2+x_0^2+y_0^2+z_0^2$ and hence $3n+15=(2w_0)^2+(2x_0)^2+(2y_0)^2+(2z_0)^2$.
So, there are nonzero integers $w,x,y,z$ with $2\mid wxyz$ such that $3n+15=w^2+x^2+y^2+z^2$.

If three of $w,x,y,z$, say $x,y,z$, are even, then two of them, say $y$ and $z$, are congruent modulo $4$.
If two of $w,x,y,z$, say $y$ and $z$, are odd, then $y\eq \ve z\pmod 4$ for a suitable choice of $\ve\in\{\pm1\}$.
So, without loss of generality, we may assume that $2\mid wx$ and $y\eq z\pmod 4$. Since $(y-z)/2\eq0\pmod 2$ and
$$y^2+z^2=2\l(\f{y+z}2\r)^2+2\l(\f{y-z}2\r)^2,\tag4.4$$
we have $3n+15=w^2+x^2+2u^2+2(2v)^2$ for some integers $u$ and $v$ not all zero.

{\it Case} 1. $3\nmid wx$.

In this case, $u^2+(2v)^2+1\eq0\pmod 3$, hence $3\nmid uv$ and we are done.

{\it Case} 2. $3\mid w$ and $3\mid x$.

In this case, we have $u\eq v\eq0\pmod 3$. By Lemma 4.2, we may write $w^2+2u^2$ as $q^2+2r^2$ with $q,r\in\Z$ and $3\nmid qr$,
and write $x^2+8v^2$ as $s^2+8t^2$ with $s,t\in\Z$ and $3\nmid st$. As $2\mid wx$, we have $2\mid qs$. Note that
$$3n+15=q^2+s^2+2r^2+2(2t)^2\ \ \t{with}\ 3\nmid qsrt.$$

{\it Case}\ 3. Exactly one of $w$ and $x$ is a multiple of $3$.

Without loss of generality, we assume that $3\nmid w$ and $3\mid x$.
Clearly $u^2+(2v)^2\eq1\pmod 3$, hence exactly one of $u$ and $v$ is a multiple of $3$.
If $3\mid u$ and $3\nmid v$, then by Lemma 4.2 we can write $x^2+2u^2$ as $r^2+2s^2$ with $r,s\in\Z$ and $3\nmid rs$.
If $3\nmid u$ and $3\mid v$, then by Lemma 4.2 we can write $x^2+8v^2$ as $r^2+8t^2$ with $r,t\in\Z$ and $3\nmid rt$.
As $r\eq x\pmod 2$ and $2\mid wx$, we have $2\mid wr$. Anyway, $3n+15=w^2+r^2+2s^2+2(2t)^2$ for some $s,t\in\Z$ with $3\nmid st$.

In view of the above,  $p_8+2p_8+4p_8+8p_8$ is indeed universal over $\Z$. It follows that
$$p_8+p_8+2p_8+8p_8,\ p_8+p_8+2p_8+2p_8,\ p_8+2p_8+2p_8+4p_8$$
are also universal over $\Z$ since $4p_8(x)+1=p_8(1-2x)$ by (1.3).

(b) Fix $d\in\{1,3,5,7,9,11,13\}$. Now we turn to show the universality of $p_8+2p_8+4p_8+dp_8$ over $\Z$
which implies the universality of $p_8+p_8+2p_8+dp_8$ over $\Z$. It suffices to show that for any $n\in\N$
we have $3n+d+7=x^2+y^2+2z^2+dw^2$ for some $w,x,y,z\in\Z$ with $2\mid xy$ and $3\nmid wxyz$.

If $n\ls d-2$ (i.e., $3n+d+7\ls4d+1$), then we may check via computer that $3n+d+7$ can be indeed written as $x^2+y^2+2z^2+dw^2$
with $w,x,y,z\in\Z$, $2\mid xy$ and $3\nmid wxyz$. For example, for $d=3$ we have
$$3\times0+10=1^2+2^2+2\times1^2+3\times1^2\ \t{and}\ 3\times1+10=2^2+2^2+2\times1^2+3\times1^2.$$

Now let $n>d-2$. Choose $w\in\{1,2\}$ with $w\not\eq n\pmod 2$. Then $3n+d+7-dw^2>1$ and $3n+d+7-dw^2\eq1\pmod 6$.
By Lemma 4.3(ii), $3n+d+7-dw^2=x^2+y^2+2z^2$ for some $x,y,z\in\Z$ with $2\mid xy$ and $3\nmid xyz$. So we have the desired result.

(c) Now we prove the universality of $p_8+2p_8+2p_8+2p_8$ over $\Z$.
It suffices to show that for any $n\in\N$ we can write $3n+7=w^2+2(x^2+y^2+z^2)$
with $w,x,y,z\in\Z$ and $3\nmid wxyz$.

Clearly,
$$7=1^2+2(1^2+1^2+1^2)\ \ \t{and}\ 10=2^2+2(1^2+1^1+1^2).$$
Now we assume $n\gs 2$. Since $p_8+p_8+2p_8+9p_8$ is universal over $\Z$, we can write
$$3n+7=3(n-2)+13=2w^2+x^2+y^2+z^2,$$
where $w,x,y,z\in\Z$, one of $x,y,z$ is divisible by $3$ but not divisible by $9$, and the other three of $w,x,y,z$ are all coprime to $3$.
Clearly, two of $x,y,z$, say $y$ and $z$, have the same parity.
As $y$ and $z$ are not all divisible by $3$, $(y+z)/2$ and $(y-z)/2$ are not all divisible by $3$.
So, in view of (4.4), there are $u,v\in\Z$ not all divisible by $3$ such that
$$3n+7=2w^2+x^2+2u^2+2v^2.$$
Without loss of generality, we suppose that $3\nmid v$. Note that $x\not=0$ and
$$x^2+2u^2\eq 7-2v^2-2w^2\eq0\pmod 3.$$
By Lemma 4.2, we can write $x^2+2u^2$ as $s^2+2t^2$ with $s,t\in\Z$ and $3\nmid st$.
Therefore,
$$3n+7=s^2+2t^2+2v^2+2w^2\quad\t{with}\ 3\nmid stvw.$$

(d) Now we prove that $p_8+2p_8+4p_8+4p_8$ is universal over $\Z$.
It suffices to show that for any $n\in\N$ we have
$3n+11=2w^2+x^2+y^2+4z^2$ for some $w,x,y,z\in\Z$ with $2\mid xy$ and $3\nmid wxyz$.

As $p_8+p_8+2p_8+4p_8$ is universal over $\Z$, there are $w,x,y,z\in\Z$ with $3\nmid wxyz$
such that $3(n+1)+8=2w^2+x^2+y^2+4z^2$. We are done if $x$ or $y$ is even.
Now assume $2\nmid xy$. Note that $(x+y)/2\not\eq(x-y)/2\pmod 2$.
Without loss of generality we suppose that $w\eq(x+y)/2\pmod 2$ since $(-y)^2=y^2$.
Clearly,
$$\align 3n+11=&2w^2+2\l(\f{x+y}2\r)^2+2\l(\f{x-y}2\r)^2+4z^2
\\=&\l(w+\f{x+y}2\r)^2+\l(w-\f{x+y}2\r)^2+2\l(\f{x-y}2\r)^2+4z^2
\endalign$$
with $w\pm(x+y)/2$ even. If $x\not\eq y\pmod 3$, then
$$\l(w+\f{x+y}2\r)^2+\l(w-\f{x+y}2\r)^2\eq11-2\l(\f{x-y}2\r)^2-4z^2\eq 11-6\eq 2\pmod 3$$
and hence $w\pm(x+y)/2\not\eq0\pmod 3$. When $x\eq y\pmod 3$, exactly one of $w+(x+y)/2$ and $w-(x+y)/2$,
is divisible by 3, and we may simply assume that $w\eq (x+y)/2\eq x\pmod 3$ (otherwise we may use $-w$ to replace $w$),
hence either $w=x=y$ or
$$\l(w-\f{x+y}2\r)^2+2\l(\f{x-y}2\r)^2=u^2+2v^2\ \t{for some}\ u,v\in\Z\ \t{with}\ 3\nmid uv$$
(in view of Lemma 4.2). If $w=x=y$, then
$$3n+11=2w^2+x^2+y^2+4z^2=4w^2+4z^2\eq0\pmod 4.$$

Now we suppose that $3n+11\eq0\pmod 4$. If $3n+11-4\times1^2=4x^2$ and $3n+11-4\times5^2=4y^2$ for some $x,y\in\N$, then
$$(x+y)(x-y)=5^2-1^2=24$$
and hence $(x,y)\in\{(5,1),(7,5)\}$, therefore
$$3n+11\in\{4(5^2+1^2),4(7^2+1^2)\}=\{104,200\}.$$
Observe that
$$104=2\times 2^2+4^2+4^2+4\times 4^2\ \ \t{and}\ \ 200=2\times 2^2+8^2+8^2+4\times 4^2.$$
Now assume that $3n+11\not=104,200$. Then $3n+11-4w^2$ is not a square for a suitable choice of $w\in\{1,5\}$.
If $3n+11\ls w^2+1$, then $3n+11=20$. Note that
$$20=2\times 2^2+2^2+2^2+4\times1^2.$$
As $3n+11-w^2\eq1\pmod6$, if $3n+11>w^2+1$ then by Lemma 4.3(ii) we have $3n+11-w^2=x^2+2y^2+4z^2$ for some $x,y,z\in\Z$ with $3\nmid xyz$,
also we don't have $|w|=|x|=|y|$ since $3n+11-4w^2$ is not a square.
Therefore we get the desired result by the arguments in the last paragraph.

(e) As $4p_8(x)+1=p_8(1-2x)$, the universality of $p_8+2p_8+2p_8+3p_8$ over $\Z$ follows from the universality of $p_8+2p_8+8p_8+3p_8$ over $\Z$.
Now we prove that $p_8+p_8+p_8+3p_8$ and $p_8+2p_8+3p_8+8p_8$ are universal over $\Z$.
Let $n$ be any nonnegative integer. It suffices to show that $3n+6=3w^2+x^2+y^2+z^2$ for some $w,x,y,z\in\Z$ with $3\nmid wxyz$,
and that $3n+14=3w^2+x^2+2y^2+8z^2$ for some $w,x,y,z\in\Z$ with $3\nmid wxyz$.

For $n=0,1,\ldots,14$, via a computer we can verify that $3n+6=3w^2+x^2+y^2+z^2$ for some $w,x,y,z\in\Z$ with $3\nmid wxyz$.
Now we simply let $n\gs 15$. By Lemma 4.4(ii),
 for a suitable choice of $w\in\{1,2,4\}$, we have $3n+6-3w^2\not\in\{4^k(8l+7):\ k,l\in\N\}$ and hence $3n+6-3w^2=x^2+y^2+z^2$ for some $x,y,z\in\Z$ (by (4.1)).
Note that $x^2+y^2+z^2=3n+6-3w^2\gs 3(n+2-4^2)>0$. By Lemma 2.2(ii), there are $\bar x,\bar y,\bar z\in\Z$ with $3\nmid \bar x\bar y\bar z$ such that
$x^2+y^2+z^2=\bar x^2+\bar y^2+\bar z^2$. Thus we have the desired representation $3n+6=3w^2+\bar x^2+\bar y^2+\bar z^2$ with $3\nmid w\bar x\bar y\bar z$.

For $n=0,1,\ldots,95$, via a computer we can verify that $3n+14=3w^2+x^2+2y^2+8z^2$ for some $w,x,y,z\in\Z$ with $3\nmid wxyz$.
Now we simply let $n\gs 96$. By Lemma 4.4(ii), there are $w_1\in\{1,2,4\}$ and $w_2\in\{5,8,10\}$ such that
$3n+14-3w_i^2\not\in\{8q-3:\ q\in\Z\}\cup\{4^k(8l+7):\ k,l\in\N\}$ for $i=1,2$. If $3n+14-3w_1^2=2x_1^2$ and $3n+14-3w_2^2=2x_2^2$ with $x_1,x_2\in\N$, then
$3\nmid x_1x_2$ and $(x_1+x_2)(x_1-x_2)=3(w_2^2-w_1^2)/2$, hence $(w_1,w_2,x_1,x_2)$ is among
$$(1,5,10,8),(2,10,13,5),(2,10,20,16),(2,10,37,35),(4,8,11,7),(4,8,19,17)$$
and thus $3n+14=3w_1^2+2x_1^2$ is among $350,770,812,2750$. (Note that $n\gs96$.)
It is easy to check that
$$\align 350=&3\times1^2+1^2+2\times13^2+8\times 1^2,\ 770=3\times 2^2+2^2+2\times 11^2+8\times8^2,
\\812=&3\times1^2+1^2+2\times 2^2+8\times 10^2,\ 2750=3\times1^2+1^2+2\times 37^2+8\times1^2.
\endalign$$
Now suppose that $3n+14\not\in\{350,770,812,2750\}$. Then, for a suitable choice of $w\in\{w_1,w_2\}$ we have $3n+14-3w^2\not=2x^2$ for any $x\in\Z$.
Note that $3n+14-3w^2>3(n+4-10^2)\gs0$.
By (4.1), there are $x,y,z\in\Z$ such that $3n+14-3w^2=x^2+y^2+z^2$.
If $x,y,z$ are all even, then two of them are congruent modulo $4$.
If two of $x,y,z$, say $y$ and $z$, are odd, then we may assume that $y\eq z\pmod 4$ since $(-z)^2=z^2$.
If exactly one of $x,y,z$ is odd and the other two even numbers are not congruent modulo $4$, then
$x^2+y^2+z^2\eq 1+0+4\eq-3\pmod 8$. As our choice of $w$ guarantees that  $3n+14-3w^2\not\eq-3\pmod 8$, we may assume that $y\eq z\pmod 4$.
Let $u=(y+z)/2$ and $v=(y-z)/4$. Then
$$3n+14-3w^2=x^2+y^2+z^2=x^2+2u^2+2(2v)^2.$$
Clearly, $u$ or $v$ is not divisible by $3$. If $3\nmid v$, then $x^2+2u^2>0$ is a multiple of $3$ and hence by Lemma 4.2 we can write $x^2+2u^2$ as $r^2+2s^2$ with $r,s\in\Z$ and $3\nmid rs$. If $3\nmid u$, then $x^2+8v^2>0$ is divisible by $3$ and hence by Lemma 4.2 we can write $x^2+8v^2$ as $r^2+8t^2$ with $r,t\in\Z$ and $3\nmid rt$.
Anyway, $3n+14=3w^2+r^2+2s^2+8t^2$ for some $r,s,t\in\Z$ with $3\nmid wrst$.

(f) Now we show the universality of $p_8+2p_8+2p_8+6p_8$ over $\Z$.
It suffices to prove that for any $n\in\N$ we can write $3n+11$ as $6w^2+x^2+2y^2+2z^2$ with $w,x,y,z\in\Z$ and $3\nmid wxyz$.
This can be easily verified via a computer for $n=0,\ldots,47$. So we simply let $n>47$.
As $6\times (2^2-1^2)\eq2\pmod 4$, we have $\{3n+11-6\times 1^2,3n+11-6\times 2^2\}\not\se\{4^k(8l+7):\ k,l\in\N\}$.
So, there is a number $w_1\in\{1,2\}$ such that $3n+11-6w_1^2\not\in\{4^k(8l+7):\,k,l\in\N\}$. Similarly, there is a number $w_2\in\{4,5\}$ such that $3n+11-6w_2^2\not\in\{4^k(8l+7):\,k,l\in\N\}$. If $3n+11-6w_1^2=2x_1^2$ and $3n+11-6w_2^2=2x_2^2$ with $x_1,x_2\in\N$, then $3\nmid x_1x_2$ and
$(x_1+x_2)(x_1-x_2)=3(w_2^2-w_1^2)$, hence
$(w_1,w_2,x_1,x_2)$ is among
$$(1,4,7,2),(1,4,23,22),(1,5,11,7),(1,5,19,17),(2,4,10,8),(2,5,8,1),(2,5,32,31)$$
and $3n+11=6w_1^2+2x_1^2$ is among $224,248,728,1064,2072$. (Note that $3n+11>3\times 47+11=152$.)
Clearly,
$$\align 224=&6\times1^2+4^2+2\times1^2+2\times10^2,\ 248=6\times1^2+4^2+2\times7^2+2\times8^2,
\\728=&6\times1^2+4^2+2\times8^2+2\times17^2,\ 1064=6\times1^2+4^2+2\times11^2+2\times20^2,
\\2072=&6\times1^2+4^2+2\times20^2+2\times25^2.
\endalign$$
Now assume that $3n+11\not=224,248,728,1064,2072$. Then, there is a number $w\in\{w_1,w_2\}$ such that $3n+11-6w^2\not=2x^2$ for any $x\in\Z$.
Note that $3n+11-6w^2>3(n+3-2\times 5^2)>0$ and $3n+11-6w^2\not= 4^k(8l+7)$ for any $k,l\in\N$. By (4.1),
there are $x,y,z\in\Z$ such that $3n+11-6w^2=x^2+y^2+z^2$. Without loss of generality, we assume that $y\eq z\pmod 2$.
In view of (4.4), there are $u,v\in\Z$ such that $3n+11-6w^2=x^2+2u^2+2v^2$. Clearly, $u$ or $v$ is not divisible by $3$.
Without loss of generality, we suppose that $3\nmid v$. Note that $x^2+2u^2>0$ is a multiple of $3$. By Lemma 4.2 we can write $x^2+2u^2$
as $s^2+2t^2$ with $s,t\in\Z$ and $3\nmid st$. Therefore $3n+11=6w^2+s^2+2t^2+2v^2$ with $3\nmid stvw$.

(g) Now we prove that $p_8+2p_8+4p_8+12p_8$ is universal over $\Z$ (which implies the universality
of $p_8+p_8+2p_8+12p_8$ over $\Z$ by (1.3)). It suffices to show that for any $n\in\N$
we can write $3n+19$ as $12w^2+x^2+2y^2+4z^2$ with $w,x,y,z\in\Z$ and $3\nmid wxyz$. By [D39, pp.\,112-113],
$$\{x^2+2y^2+4z^2:\ x,y,z\in\Z\}=\N\sm\{4^k(16l+14):\ k,l\in\N\}.\tag4.5$$
As $12(2^2-1^2)\eq12(5^2-4^2)\eq4\pmod 8$, there are $w_1\in\{1,2\}$ and $w_2\in\{4,5\}$ such that $3n+19-12w_i^2\not\in\{4^k(16l+14):\,k,l\in\N\}$ for $i=1,2$.
If $3n+19-12w_1^2=x_1^2$ and $3n+19-12w_2^2=x_2^2$ with $x_1,x_2\in\N$, then $3\nmid x_1x_2$ and
$(x_1+x_2)(x_1-x_2)=12(w_2^2-w_1^2)$, hence
$(w_1,w_2,x_1,x_2)$ is among
$$\align&(1,4,14,4),(1,4,46,44),(1,5,17,1),(1,5,22,14),(1,5,38,34),(1,5,73,71),
\\&(2,4,13,5),(2,4,20,16),(2,4,37,35),(2,5,16,2),(2,5,64,62)
\endalign$$
and $3n+19=12w_1^2+x_1^2$ belongs to the set
$$T=\{208,\ 217,\ 301,\ 304,\ 448,\ 496,\ 1417,\ 1456,\ 2128,\ 4144,\ 5341\}.$$
Via a computer we can check that each element of $T\cup\{3n+19:\ n=0,\ldots,93\}$ can be written as $12w^2+x^2+2y^2+4z^2$ with $w,x,y,z\in\Z$ and $3\nmid wxyz$.

Now we may suppose that $n>93$ and that there is a number $w\in\{w_1,w_2\}$ such that $3n+19-12w^2$ is not a square.
Note that $3n+19-12w^2>3(n+6-4\times 5^2)\gs0$ and $3n+19-12w^2\not=4^k(16l+14)$ for any $k,l\in\N$. By (4.5), there are $x,y,z\in\Z$ with $2\mid xy$
such that $3n+19-12w^2=x^2+y^2+2z^2$. Clearly, $x$ or $y$ is not divisible by $3$. Without loss of generality, we suppose that $3\nmid x$.
Note that $y^2+2z^2>0$ is a multiple of $3$. By Lemma 4.2, we may write $y^2+2z^2=u^2+2v^2$ with $u,v\in\Z$ and $3\nmid uv$. Thus
$3n+19=12w^2+x^2+u^2+2v^2$ with $3\nmid uvwx$ and $ux\eq xy\eq0\pmod 2$.

(h) Let $d\in\{3,5\}$. Now we prove that $p_8+2p_8+4p_8+2dp_8$ is universal over $\Z$
(which implies that $p_8+p_8+2p_8+2dp_8$ is also universal over $\Z$). It suffices to show that
for any $n\in\N$ we can write $3n+2d+7$ as $2dw^2+x^2+2y^2+4z^2$ with $w,x,y,z\in\Z$ and $3\nmid wxyz$.
As $7\eq 3\eq-5\pmod 4$ and $7\not\eq -5\pmod 8$, by Lemma 4.4 there are $w_1\in\{1,2,4\}$ and $w_2\in\{5,8,10\}$
such that for $i=1,2$ we have $(3n+2d+7)/2-dw_i^2\not\in\{4^k(8l+7):\, k,l\in\N\}$ and hence
$$3n+2d+7-2dw_i^2\not\in\{4^k(16l+14):\ k,l\in\N\}.$$
If $3n+2d+7-2dw_1^2=x_1^2$ and $3n+2d+7-2dw_2^2=x_2^2$ with $x_1,x_2\in\N$, then $3\nmid x_1x_2$
and $(x_1+x_2)(x_1-x_2)=2d(w_2^2-w_1^2)$,
hence $(d,w_1,w_2,x_1,x_2)$ is among
$$\align&(3,1,5,13,5),(3,1,5,20,16),(3,1,5,37,35),(3,2,8,19,1),(3,2,8,23,13),
\\&(3,2,8,47,43),(3,2,8,91,89),(3,2,10,25,7),(3,2,10,26,10),(3,2,10,40,32),
\\&(3,2,10,74,70),(3,2,10,145,143),(3,4,8,17,1),(3,4,8,22,14),(3,4,8,38,34),
\\&(3,4,8,73,71),(3,4,10,23,5),(3,4,10,25,11),(3,4,10,65,61),(3,4,10,127,125);
\\&(5,1,5,16,4),(5,1,5,17,7),(5,1,5,19,11),(5,1,5,23,17),(5,1,5,32,28),
\\&(5,1,5,61,59),(5,2,8,25,5),(5,2,8,31,19),(5,2,8,35,25),(5,2,8,53,47),
\\&(5,2,8,77,73),(5,2,8,151,149),(5,2,10,31,1),(5,2,10,32,8),(5,2,10,34,14),
\\&(5,2,10,38,22),(5,2,10,46,34),(5,2,10,53,43),(5,2,10,64,56),(5,2,10,83,77),
\\&(5,2,10,122,118),(5,2,10,241,239),(5,4,8,22,2),(5,4,8,23,7),(5,4,8,26,14),
\\&(5,4,8,34,26),(5,4,8,43,37),(5,4,8,62,58),(5,4,8,121,119),(5,4,10,29,1),
\\&(5,4,10,31,11),(5,4,10,37,23),(5,4,10,41,29),(5,4,10,73,67),(5,4,10,107,103),
\\&(5,4,10,211,209)
\endalign$$
and thus $3n+2d+7=2dw_1^2+x_1^2\in E(d)$, where
$$\align E(3)=\{&175, 385, 406, 553, 580, 625, 649, 700, 721, 1375, 1540,
\\&1624, 2233, 4321, 5425, 5500, 8305,16225, 21049\}
\endalign$$
and
$$\align E(5)=\{&266, 299, 371, 539, 644, 665, 689, 836, 1001, 1034, 1064,1121,
\\& 1196, 1265, 1316, 1484, 1529,1841, 2009, 2156, 2849, 3731, 4004,
\\&4136, 5489, 5969, 6929, 11609, 14801, 14924, 22841, 44681,58121\}.
\endalign$$
If $3n+2d+7<200d$ or $3n+2d+7\in E(d)$, then we may use a computer to check that
$3n+2d+7$ can be indeed written as $2dw^2+x^2+2y^2+4z^2$ with $w,x,y,z\in\Z$ and $3\nmid wxyz$.
For example,
$$21049=6\times1^2+7^2+2\times53^2+4\times 62^2,\ 58121=10\times2^2+65^2+2\times100^2+4\times92^2.$$

 Now let $3n+2d+7\gs 200d$ and $3n+2d+7\not\in E(d)$. Then, for some $w\in\{w_1,w_2\}$ the number
 $3n+2d+7-2dw^2$ is not a square. Clearly, $3n+2d+7-2dw^2\gs 2d(10^2-w^2)\gs0$. As $3n+2d+7-2dw^2
 \not=4^k(16l+14)$ for any $k,l\in\N$, by (4.5) there are $x,y,z\in\Z$ with $2\mid xy$ such that $3n+2d+7-2dw^2=x^2+y^2+2z^2$.
As $x^2+y^2+2z^2\eq 7\eq1\pmod 3$, $x$ or $y$ is not divisible by $3$. Without loss of generality, we assume that $3\nmid x$.
Note that $y^2+2z^2>0$ is a multiple of $3$. By Lemma 4.2, we can write $y^2+2z^2$ as $u^2+2v^2$ with $u,v\in\Z$ and $3\nmid uv$.
Thus $3n+2d+7=2dw^2+x^2+u^2+2v^2$ with $3\nmid uvwx$ and $ux\eq xy\eq0\pmod 2$.

(i) Let $d\in\{2,3\}$. We now show that $dp_8+p_8+2p_8+5p_8$ is universal over $\Z$.
It suffices to prove that for any $n\in\N$ we can write $3n+d+8$ as $dw^2+x^2+2y^2+5z^2$ with $w,x,y,z\in\Z$ and $3\nmid wxyz$.
As $1^2,2^2,5^2$ are pairwise incongruent modulo $5$, there are $w_1,w_2\in\{1,2,5\}$ with $w_1<w_2$ such that $3n+d+8-dw_i^2\not\eq0\pmod 5$
for $i=1,2$.

Suppose that $3n+d+8-dw_1^2=2x_1^2$ and $3n+d+8-dw_2^2=2x_2^2$ with $x_1,x_2\in\N$. Then
$2(x_1+x_2)(x_1-x_2)=d(w_2^2-w_1^2)$. In the case $(d,w_1,w_2)=(2,1,2)$, we have $(x_1,x_2)=(2,1)$ and hence
$3n+d+8=dw_1^2+2x_1^2=2+2\times2^2=10$. If $(d,w_1,w_2)=(2,1,5)$, then $(x_1,x_2)$ is $(5,1)$ or $(7,5)$, hence
$3n+d+8=dw_1^2+2x_1^2=2(x_1^2+1)\in\{52,100\}$. If $(d,w_1,w_2)=(2,2,5)$, then $(x_1,x_2)$ is $(5,2)$ or $(11,10)$, hence
$3n+d+8=dw_1^2+2x_1^2=2(2^2+x_1^2)\in\{58,250\}$. When $d=3$, we must have $2\mid(w_1^2-w_2^2)$ and $3\nmid x_1x_2$, hence $(w_1,w_2)=(1,5)$ and $(x_1,x_2)=(10,8)$, therefore
$3n+d+8=dw_1^2+2x_1^2=3+2\times10^2=203$.
Note that
$$\align 10=&2\times 1^2+1^2+2\times1^2+5\times 1^2,\ 52=2\times 2^2+4^2+2\times2^2+5\times 2^2,
\\58=&2\times1^2+1^2+2\times 5^2+5\times1^2,\ 100=2\times1^2+4^2+2\times 1^2+5\times4^2,
\\250=&2\times 1^2+1^2+2\times1^2+5\times7^2,\ 203=3\times 1^2+5^2+2\times5^2+5\times 5^2.
\endalign$$

Now we handle the remaining case. Assume that $(3n+d+8-dw^2)/2$ is not a square for a suitable choice of $w\in\{w_1,w_2\}$.
If $3n+d+8<25d$, then $n\ls 13$ for $d=2$, and $n\ls 21$ for $d=3$. Thus, when $3n+d+8<25d$, via a computer we can write
$3n+d+8$ in the form $dv^2+x^2+2y^2+5z^2$ with $v,x,y,z\in\Z$ and $3\nmid vxyz$. Now let $3n+d+8\gs 25d\gs dw^2$.
As $3n+d+8-dw^2\not\eq0\pmod 5$ and
$$\{x^2+2y^2+5z^2:\ x,y,z\in\Z\}=\N\setminus\{25^k(25l\pm 10):\ k,l\in\N\}\tag4.6$$
by [D39, pp. 112-113], there are $x,y,z\in\Z$ such that $3n+d+8-dw^2=x^2+2y^2+5z^2$.
Note that $3n+d+8-dw^2\eq2\pmod 3$ and $3n+d+8-dw^2\not=2y^2$. If $3\nmid y$, then $3\mid (x^2+5z^2)$,
hence by Lemma 4.2 we can write $x^2+5z^2=u^2+5v^2$ with $u,v\in\Z$ and $3\nmid uv$. In the case $3\mid y$, we must have
$3\mid x$ and $3\nmid z$, hence we may write $x^2+2y^2>0$ as $u^2+2v^2$ with $u,v\in\Z$ and $3\nmid uv$. Anyway, $3n+d+8
=dw^2+r^2+2s^2+5t^2$ for some $r,s,t\in\Z$ with $3\nmid rst$.

(j) To prove the universality of $p_8+p_8+3p_8+5p_8$ over $\Z$, it suffices to show that for any $n\in\N$ we have
$3n+10-3w^2=x^2+y^2+5z^2$ for some $w,x,y,z\in\Z$ with $3\nmid wxyz$.

It is known (cf. [D39, pp. 112-113]) that
$$\{x^2+y^2+5z^2:\ x,y,z\in\Z\}=\N\sm\{4^k(8l+3):\ k,l\in\N\}.\tag4.7$$
By Lemma 4.4, there are $w_1\in\{1,2,4\}$ and $w_2\in\{5,8,10\}$ such that
$$\{3n+10-3w_1^2,\ 3n+10-3w_2^2\}\cap\{4^k(8l+3):\ k,l\in\N\}=\emptyset.$$
If $3n+10-3w_1^2=x_1^2$ and $3n+10-3w_2^2=x_2^2$ with $x_1,x_2\in\N$, then $3\nmid x_1x_2$ and
$$(x_1+x_2)(x_1-x_2)/3=w_2^2-w_1^2\in\{y^2-x^2: x\in\{1,2,4\}\ \&\ y\in\{5,8,10\}\},$$
hence $(w_1,w_2,x_1,x_2)$ is among
$$\align&(1,5,11,7),(1,5,19,17),(2,5,8,1),(2,5,32,31),(4,5,14,13),
\\&(1,8,17,10),(1,8,95,94),(2,8,14,4),(2,8,46,44),(4,8,13,5),
\\&(4,8,20,16),(4,8,37,35),(1,10,19,8),(1,10,149,148),(2,10,17,1),
\\&(2,10,22,14),(2,10,38,34),(2,10,73,71),(4,10,16,2),(4,10,64,62)
\endalign$$
and it follows that $3n+10=3w_1^2+x_1^2$ belongs to the set
$$\align E=&\{76, 124, 208, 217, 244, 292, 301, 304, 364, 448, 496,
\\&\ 1036, 1417, 1456, 2128, 4144, 5341,9028, 22204\}.
\endalign$$
Via a computer we can write each element of $E\cup\{3n+10:\ n=0,\ldots,96\}$ in the form $3w^2+x^2+y^2+5z^2$ with $w,x,y,z\in\Z$ and $3\nmid wxyz$;
for example,
$$76=3\times1^2+2^2+7^2+5\times2^2\ \t{and}\ 22204=3\times1^2+20^2+26^2+5\times 65^2.$$

Now suppose that $n\gs 97$ and  $3n+10\not\in E$. Then there is a suitable choice of $w\in\{w_1,w_2\}$ such that
$3n+10-3w^2$ is not a square. Clearly, $3n+10=3(n+3)+1>3\times 10^2\gs 3w^2$ and $3n+10-3w^2\not=4^k(8l+3)$ for any $k,l\in\N$.
Thus, by (4.7) there are $x,y,z\in\Z$ such that $3n+10-3w^2=x^2+y^2+5z^2$. Clearly $x$ or $y$ is not divisible by $3$.
Without loss of generality, we assume that $3\nmid x$. Then $y^2+5z^2$ is a positive integer divisible by $3$.
Applying Lemma 4.2 we find that $y^2+5z^2=u^2+5v^2$ for some $u,v\in\Z$ with $3\nmid uv$. Therefore
$3n+10=3w^2+x^2+u^2+5v^2$ with $3\nmid uvwx$. This concludes the proof. \qed

\heading{5. Conjectures involving $p_m(x)$ with $m\in\{5,6,7\}$}\endheading

In 2008 the author (cf. [S09]) conjectured that $216$ is the only natural number which cannot be written as $p+x(x+1)/2$,
where $p$ is prime or zero, and $x$ is an integer. Here we pose a similar conjecture involving generalized pentagonal numbers.

\proclaim{Conjecture 5.1} Any $n\in\N$ can be expressed as $p+x(3x-1)/2$, where $p$ is an odd prime or zero, and $x$ is an integer.
In other words, any nonnegative integer is either an odd prime, or a generalized pentagonal number, or the sum of an odd prime and a generalized pentagonal number.
\endproclaim
\Remark\ 5.1. We have verified Conjecture 5.1 for all $n=0,\ldots,10^9$.
\medskip

For any $m\in\{5,6,7,\ldots\}$, we define $\bar p_m(x):=p_m(-x)$. Those
$$\bar p_m(n)=p_m(-n)=(m-2)\f{n(n+1)}2-n\ \ (n=0,1,2,\ldots)$$
are usually called the {\it second $m$-gonal numbers}.

\proclaim{Conjecture 5.2} {\rm (i)} Both $p_5+p_5+\bar p_5+\bar p_5$ and $p_5+p_5+p_5+\bar p_5$
are universal over $\N$. Moreover, for any $n\in\Z^+$, there are $w,x,y,z\in\N$ with $x,y,z$ not all even such that
$n=\bar p_5(w)+p_5(x)+p_5(y)+ p_5(z).$

{\rm (ii)} $p_5+bp_5+cp_5+dp_5$ is universal over $\N$ if $(b,c,d)$ is among the following $15$ triples:
$$\gather(1,1,2),(1,2,2),(1,2,3),(1,2,4),(1,2,5),(1,2,6),(1,3,6),
\\(2,2,4),(2,2,6),(2,3,4),(2,3,5),(2,3,7),(2,4,6),(2,4,7),(2,4,8).
\endgather$$
\endproclaim

\proclaim{Conjecture 5.3}  All the sums
$$\align &p_6+p_6+\bar p_6+\bar p_6,\ p_6+p_6+2p_6+4p_6,\ p_6+2p_6+\bar p_6+\bar p_6,
\\&p_6+2p_6+\bar p_6+2\bar p_6,\ p_6+p_6+2p_6+\bar p_6,\ p_6+p_6+3p_6+\bar p_6,
\\&p_6+p_6+4p_6+\bar p_6,\ p_6+p_6+8p_6+\bar p_6,\ p_6+2p_6+2p_6+\bar p_6,
\\&p_6+2p_6+3p_6+\bar p_6,\ p_6+2p_6+3p_6+2\bar p_6,\ p_6+2p_6+4p_6+\bar p_6
\endalign$$
are universal over $\N$.
\endproclaim

\proclaim{Conjecture 5.4} {\rm (i)} Any $n\in\N$ with $n\not=23$ can be written as $p_7(x)+p_7(y)+2p_7(z)$ with $x,y,z\in\Z$.
Also, all the sums
$$\align &p_7+\bar p_7+2p_7+2\bar p_7,\ \bar p_7+p_7+p_7+2p_7,\ \bar p_7+p_7+p_7+3p_7,
\\& \bar p_7+p_7+2p_7+3p_7,\ \bar p_7+p_7+2p_7+8p_7
\endalign$$ are universal over $\N$.

{\rm (ii)} For each $m\in\{7,9,10,11,12,13,14\}$, the sum $p_m+2p_m+4p_m+8p_m$ is universal over $\Z$.
\endproclaim

\Remark\ 5.2. Guy [Gu] noted that none of $10,\, 16,\, 76$ can be written as the sum of three generalized heptagonal numbers.
We guess that
$$\{p_7(x)+p_7(y)+p_7(z):\ x,y,z\in\Z\}=\N\sm\{10,16,76,307\}$$
and
$$\{p_7(x)+2p_7(y)+4p_7(z):\ x,y,z\in\Z\}=\N\sm\{131,146\}.$$

\medskip

\Ack. The author thanks his graduate student X.-Z. Meng for helpful comments on Example 3.1.

\enddocument